\documentclass[a4paper]{elsarticle}
\usepackage[utf8]{inputenc}
\usepackage{import}

\usepackage{fontawesome}

\usepackage{fullpage}

\usepackage{lipsum}

\usepackage{circledsteps}
\pgfkeys{/csteps/fill color=white}
\pgfkeys{/csteps/inner ysep=5pt}
\pgfkeys{/csteps/inner xsep=5pt}

\usepackage[english]{babel}

\usepackage{enumitem}
\setlist[1]{itemsep=-5pt}
\usepackage[margin=2.59cm]{geometry}
\usepackage{amsmath,amssymb,amsthm}
\usepackage{cases}
\usepackage{mathrsfs}
\usepackage{bm}
\newtheorem{remark}{Remark}[section]
\usepackage{algorithm}
\usepackage{algpseudocode}
\theoremstyle{definition}

\usepackage{stmaryrd}
\usepackage{cases}

\usepackage{appendix}
\providecommand{\appendixname}{Appendix}
\usepackage{pdflscape}
\usepackage{tabularx}
\usepackage{booktabs}

\usepackage{multicol}
\usepackage{multirow}
\usepackage{xparse}
\usepackage{xspace}
\usepackage{makecell}

\usepackage{xcolor}
\usepackage{footnote} 
\usepackage{graphicx}
\usepackage{svg}
\usepackage{epstopdf}
\usepackage{tikz}
\usetikzlibrary{calc}
\usepackage{pgfplots}
\pgfplotsset{compat=newest}

\usepackage[percent]{overpic}
\usepackage{placeins} 

\usepackage[colorlinks=true]{hyperref}
\usepackage{cleveref}
\biboptions{sort&compress}
\usepackage[numbers]{natbib}
\usepackage{caption}

\usepackage[notref,notcite,final]{showkeys}

\usepackage[color=gray, backgroundcolor=yellow, textwidth=2cm, textsize=footnotesize]{todonotes}
\usepackage{setspace}
\usepackage{subcaption}
\usepackage{fullpage}
\usepackage{framed}
\usepackage[noprefix]{nomencl}

\makenomenclature

\renewcommand*\nompreamble{\begin{multicols}{2}}
\renewcommand*\nompostamble{\end{multicols}}
\usepackage{etoolbox}
\usepackage{mathtools}
\usepackage{pdfpages}

\newcommand{\tr}{\mathrm{tr}\,}
\newcommand{\sig}{\bm{\sigma}}

\newcommand{\squad}{\hspace{0.5em}}
\newcommand{\Bu}{\bm{u}}
\newcommand{\eps}{\bm{\varepsilon}}

\newcommand{\epsp}{\bm{\varepsilon}^\mathrm{p}}
\newcommand{\epsr}{\bm{\varepsilon}^\mathrm{r}}
\newcommand{\setI}{\mathrm{I}}

\newcommand{\discrete}[1]{\mkern 1.2mu\underline{\mkern-1.2mu#1\mkern-1.2mu}\mkern 1.2mu}

\makeatletter
\newcommand{\multiline}[1]{%
  \begin{tabularx}{\dimexpr\linewidth-\ALG@thistlm}[t]{@{}X@{}}
    #1
  \end{tabularx}
}
\makeatother

\newcommand{\tB}[1]{{\color{black} #1}}

\begin{document}

\begin{frontmatter}

\title{
Multi-temporal decomposition for elastoplastic ratcheting solids
}

\author[add1]{Jacinto Ulloa\corref{cor1}}
\ead{julloa@caltech.edu}
\author[add2]{Geert Degrande}
\ead{geert.degrande@kuleuven.be}
\author[add1]{José E. Andrade}
\ead{jandrade@caltech.edu}
\author[add2]{Stijn François\corref{cor1}}
\ead{stijn.francois@kuleuven.be}

\cortext[cor1]{Corresponding author}
\address[add1]{Division of Engineering and Applied Science, California Institute of Technology, Pasadena, CA 91125, USA}
\address[add2]{Department of Civil Engineering, KU Leuven, B-3001 Leuven, Belgium}

\begin{abstract}
This paper presents a multi-temporal formulation for simulating elastoplastic solids under cyclic loading. We leverage the proper generalized decomposition (PGD) to decompose the displacements into multiple time scales, separating the spatial and intra-cyclic dependence from the inter-cyclic variation. In contrast with the standard incremental approach, which solves the (non-linear and computationally intensive) mechanical balance equations at every time step, the proposed PGD approach allows the mechanical balance equations to be solved exclusively for the small-time intra-cyclic response, while the large-time inter-cyclic response is described by simple scalar algebraic equations. Numerical simulations exhibiting complex cyclic responses, including a 2D problem and an application to a monopile foundation, demonstrate that PGD solutions with a limited number of space-time degrees of freedom may be obtained numerically, only requiring a few modes to accurately capture the reference response.
\end{abstract}


\end{frontmatter}

\section{Introduction}
\label{sec_intro}
\FloatBarrier

Challenges in modeling non-linear solids under cyclic loading emerge from the need to describe complex behavior such as cyclic plasticity and fatigue, and from the computational expense of simulating a large number of loading cycles. 
Advanced constitutive models have been developed to describe the cyclic response of materials,
with applications including granular soils~\cite{suiker2003,niemunis2005,franccois2010,houlsby2017} and ductile solids~\cite{armstrong1966,chaboche2008,ulloa2021b}. These models are highly non-linear and are typically solved using finite element techniques for spatial discretization of the mechanical balance equation in conjunction with time integration schemes for internal variables. The resulting incremental problem generally requires small increments to achieve numerical convergence and accurately capture the non-linear response, resulting in a high computation time for problems with a large number of degrees of freedom. Under cyclic loading conditions, the computational cost scales with the number of cycles and further poses storage requirements to track the full-time history of the~response.

Integration methods that avoid a step-by-step calculation of the full-time history have been proposed in the literature for models with internal variables. Examples of non-intrusive procedures are cycle jumping techniques \citep{van2001jump, cojo2006jump} and the Fourier-transformation temporal integration method \citep{ftti1, ftti2}. These methods perform step-by-step calculations for a few cycles and consider an extrapolation in between. While the simplicity of these methods is appealing, the definition of the jump size may not be trivial and small jumps may be required for complex material behavior. A more intrusive procedure that has been widely used in non-linear solid mechanics is the large-time increment method (LATIN) proposed by \citet{lade1989}, where the global equilibrium equations are solved for the entire loading process at once. Then, the full-time history of the local constitutive model is resolved and iterated with the global solution until convergence. Applications of this framework to cyclic loading have been proposed, e.g., in the context of viscoplasticity~\cite{cognard1993} or damage~\cite{bhattacharyya2018a}. Although LATIN can handle non-linearities, its intrusive nature renders the implementation of complex material models cumbersome. 

A key aspect of LATIN is a model order reduction technique (MOR) that avoids a step-by-step solution of the global equilibrium equations. For this purpose, a modal decomposition of the state variables is performed, with each mode consisting of the product of a priori unknown spatial and temporal functions. This solution strategy can be viewed as a particular case of the proper generalized decomposition (PGD) introduced by \citet{ammar2006} for modeling complex fluids. In the PGD framework, the state variables can be separated into a priori unknown functions of space, time, and other parameters such as material properties and boundary conditions. Thus, PGD can tackle various computationally demanding problems, including multi-dimensional systems and, in general, problems that involve a large number of degrees of freedom. For instance, in close relation to the present study, references \cite{ammar2012,pasquale2021} propose a separation of the time axis into a two-dimensional space. We refer the reader to references \cite{chinesta2010,chinesta2011,chinesta2014} for further details and applications of PGD.

Exploiting the PGD formalism for non-linear models in solid mechanics is a challenging task that has not been fully addressed. Indeed, PGD formulations have been mostly limited to linear problems, while the LATIN approach is usually adopted for non-linear models. Nonetheless, a few exceptions exist in the literature. \tB{A single-scale space-time decomposition for geometrically non-linear systems was recently presented in~\citet{arjoune2022}}. On the other hand, concerning material non-linearity, \citet{bergheau2016} proposed a single-scale space-time separation of the displacement field and considered the Armstrong-Frederick~\cite{armstrong1966} plasticity model, capturing quasi-periodic behavior for a few loading cycles. A similar approach was adopted by \citet{nasri2018} for viscoplastic polycrystals and by~\citet{shirafkan2020} for elastoplastic structures. Further, a separation into a small-time scale and a large-time scale for elastoplasticity was presented in the recent work of~\citet{pasquale2023}.

In the present work, we extend the use of PGD to systems with an arbitrary number of temporal scales and complex cyclic material behavior. As a first step, we focus not on the issue of computational time but rather on finding numerical solutions admitting a multi-temporal space-time decomposition for complex systems where the response under periodic loading is not perfectly periodic. To this end, we apply the proposed formulation to a kinematic-isotropic hardening plasticity model with ratcheting, where strong cyclic-loading effects are observed.
	
\section{Multi-temporal space-time PGD}	
\label{sec:form}

\subsection{Problem setup and classical formulation}
\label{sec:setup}

Consider a solid body occupying a spatial domain $\Omega\subset\mathbb{R}^d$ and evolving in a pseudo-time interval $\setI\coloneqq[0,t_\mathrm{max}]$. The spatial boundary $\Gamma$ consists of a Dirichlet part ${\Gamma_\mathrm{D}}$ with vanishing displacements and a Neumann part $\Gamma_\mathrm{N}$ with imposed tractions $\bar{\bm{t}}(\bm{x},t)\in\mathbb{R}^d$, such that ${\Gamma_\mathrm{D} \cup \Gamma_\mathrm{N} = \Gamma}$ and ${\Gamma_\mathrm{D} \cap \Gamma_\mathrm{N} =\varnothing}$. The solid may be subjected to body forces per unit volume $\bm{b}(\bm{x},t)\in\mathbb{R}^d$. 

In the small strain setting, the response of the system is characterized by the displacement field $\Bu\colon\Omega\times\setI\to\mathbb{R}^d$, with the compatible strain tensor $\boldsymbol{\varepsilon}\colon\Omega\times\setI\to\mathbb{R}^{3\times 3}_\mathrm{sym}\coloneqq\{\bm{e}\in \mathbb{R}^{3\times 3} \ | \  \bm{e}=\bm{e}^\mathrm{T}\}$ obeying the linear relation $\boldsymbol{\varepsilon}=\nabla^\mathrm{s}\boldsymbol{u}$, where $\nabla^\mathrm{s}$ is the symmetric gradient operator. The Cauchy stress tensor $\boldsymbol{\sigma}\colon\Omega\times\setI\to\mathbb{R}^{3\times 3}_\mathrm{sym}$ is statically admissible, satisfying mechanical balance  $\forall \, t\in\setI$: 
\begin{equation}
\begin{dcases}\mathrm{div}\,\sig+\bm{b} = \boldsymbol{0} \quad &\text{in} \squad \Omega,\\
\sig\cdot\boldsymbol{n}=\bar{\bm{t}} \quad &\text{on} \squad \Gamma_{\mathrm{N}}, \\ \bm{u}=\bm{0} \quad &\text{on} \squad \Gamma_{\mathrm{D}}.\end{dcases}
\label{eq:st_adm0}
\end{equation}

The boundary value problem~\eqref{eq:st_adm0} requires a constitutive stress-strain law, which may involve dissipative inelastic behavior. Within the framework of continuum mechanics with internal variables and generalized standard materials~\cite{germain1983,halphen1975}, two energy functions characterize the system at the material level. The first is the free density $\psi\coloneqq{\psi}(\bm{\varepsilon},\mathbf{a})$, where $\mathbf{a}\colon\Omega\times\setI\to \mathbb{R}^{m}$ is a generic set of internal variables, to be specified in section~\ref{sec:const}, which account for inelastic microstructural processes phenomenologically. The second law of thermodynamics then yields the Clausius-Planck inequality 
\begin{equation}
\phi\coloneqq\sig:\dot{\bm{\varepsilon}} -\dot{\psi}(\bm{\varepsilon},\mathbf{a})\geq0.
\label{eq:CD}
\end{equation}
The Coleman-Noll procedure yields the stress-strain relation and the generalized stresses
\begin{equation}
\sig=\frac{\partial\psi(\bm{\varepsilon},\mathbf{a})}{\partial{\bm{\varepsilon}}} \quad \text{and} \quad \mathbf{s}=-\frac{\partial\psi}{\partial{\mathbf{a}}}(\bm{\varepsilon},\mathbf{a}).
\label{eq:constitutive}
\end{equation}
The second energy function follows for rate-independent systems as the functional form of the dissipation rate~\eqref{eq:CD}, i.e., $\phi=
\phi(\dot{\mathbf{a}};\mathbf{s})\geq0$. Here, $\phi$ is a convex 1-homogeneous function of $\dot{\mathbf{a}}$, non-differentiable and vanishing at the origin, and depending on the generalized stress state $\mathbf{s}$ in the case of non-associative models~\cite{ulloa2021a}. The evolution equation for $\mathbf{a}$ then follows as
\begin{equation}
\frac{\partial\psi}{\partial{\mathbf{a}}}(\bm{\varepsilon},\mathbf{a}) + \partial_{\dot{\mathbf{a}}}\phi(\dot{\mathbf{a}};\mathbf{s})\ni \bm{0},
\label{eq:biot}
\end{equation}
where $\partial_\diamond\Box(\diamond)$ is the subdifferential of a non-differentiable function $\Box$ at $\diamond$. The coupled equations~\eqref{eq:st_adm0} and~\eqref{eq:biot} represent the governing equations of the inelastic solid, to be solved for $\{\bm{u},\mathbf{a}\}$ in $\Omega\times\setI$.

The conventional approach to the problem described above starts with a discretization of the time interval $\setI$ into a finite number of time steps $0 =t_1<\dots <t_n<\dots< t_{N_t} = t_\mathrm{max}$. The solution is then sought sequentially at $t_{n}$ with all quantities known at $t_{n-1}$ for all $n\in\mathsf{I}\setminus\{1\}$, with $\mathsf{I} := \{1,\dots,N_t\}$. Using the notation $\Box_n\coloneqq\Box(\cdot\,,t_n)$, the coupled problem can be expressed as follows: $\forall\,n\in\mathsf{I}$, find $\bm{u}_n\in\mathscr{U}_{{0}}\coloneqq\{\bm{w}\in\mathrm{H}^1(\Omega;\mathbb{R}^d) \ \vert \ \bm{w}=\bm{0}\ \text{on} \ \Gamma_\mathrm{D} \}$ such that
\begin{equation}
\int_{\Omega}\sig_n:\nabla^\mathrm{s}\tilde{\bm{u}}_n\,\mathrm{d}\bm{x}-\int_{\Omega}\boldsymbol{b}_n\cdot\tilde{\bm{u}}_n\,\mathrm{d}\bm{x} -  \int_{\Gamma_\mathrm{N}}\bar{\bm{t}}_n\cdot\tilde{\bm{u}}_n\,\mathrm{d}S=0 \qquad \forall\,\tilde{\bm{u}}_n\in\mathscr{U}_{{0}},
\label{eq:st_adm1}
\end{equation}
where $\sig_n$ is given by equation~\eqref{eq:constitutive}, with $\mathbf{a}_n$ satisfying the evolution equation~\eqref{eq:biot} incrementally. The weak form~\eqref{eq:st_adm1} is non-linear due to the dependence of $\sig_n$ on $\mathbf{a}_n$; this problem is solved via Newton-Raphson schemes, which often require a fine temporal discretization for convergence. The evolution equation~\eqref{eq:biot} is solved simultaneously using standard time integration algorithms for constitutive models~\cite{simo1998,de2011,de2012}.

\subsection{Multi-temporal space-time integration}

The numerical problem described above involves solving the non-linear PDE~\eqref{eq:st_adm1} at each time step, a generally non-desirable task for long-term systems. This issue is especially pertinent for cyclic loading scenarios, where capturing the non-linear response within a single cycle may require many time steps, scaling the computational requirements drastically as the number of cycles increases. 

\begin{figure}[t]
\includegraphics[width=1\linewidth,trim={0 1cm 0 0},clip]{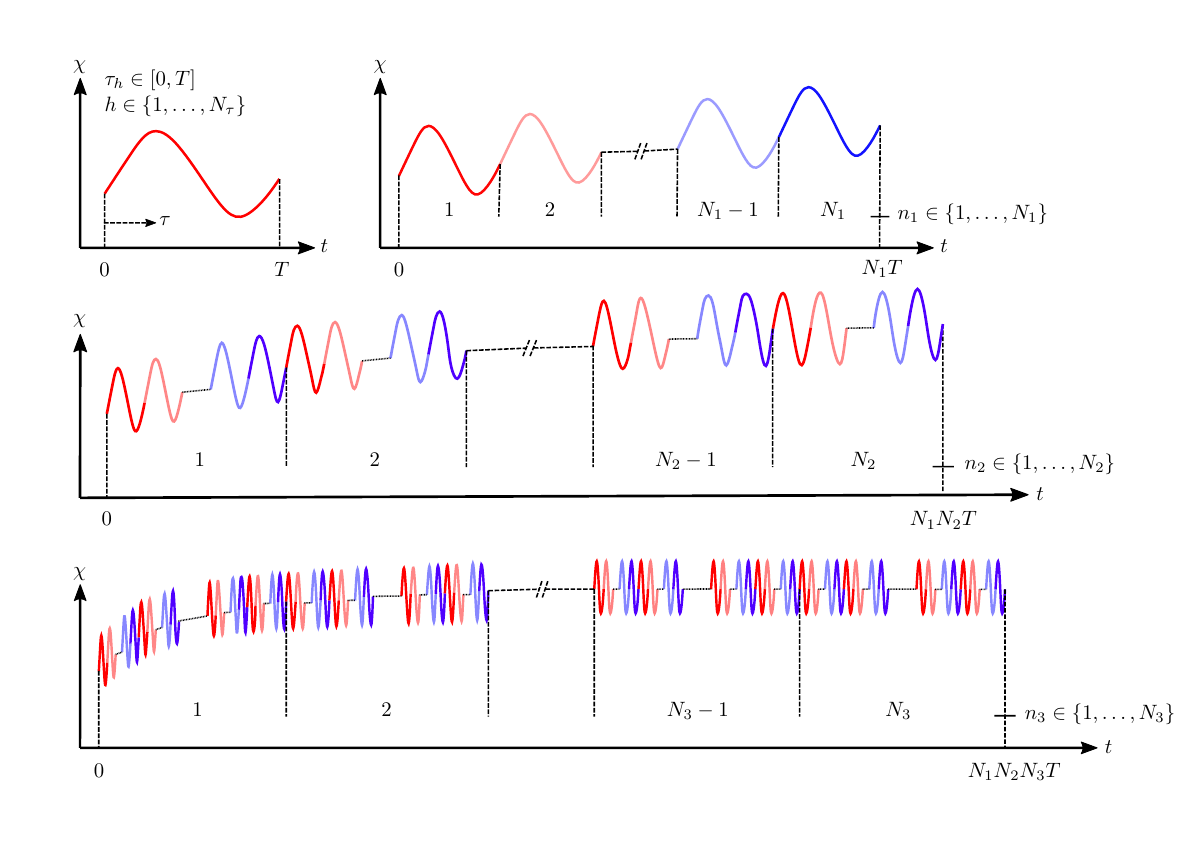}	
\caption{Example of the multi-temporal scheme using 4 time scales: a small-time scale~${\tau_h\in[0,T]}$ and 3 large-time~scales.}
\label{fig:tsep}
\end{figure}

Let us consider such a scenario, imposed by at least quasi-periodic loading conditions with period~$T$. A small-time scale is defined within each cycle at $N_\tau$ instants $0 =\tau_1<\dots <\tau_h<\dots< \tau_{N_\tau} = T$, indexed by $h\in\mathsf{T}:=\{1,\dots,N_\tau\}$, where a strongly non-linear evolution of state variables is expected.\footnote{The definition of a small-time scale with a fast evolution of internal variables in problems with cyclic loading has been used before in the context of LATIN \cite{cognard1993} and also constitutes the main idea of cycle jumping techniques \cite{van2001jump, cojo2006jump,ftti1,ftti2}.} We further consider a series of $S$ large-time scales characterized by counters ${n}_j\in\mathsf{N}_j:=\{1,\dots,N_j\}$ for all $j\in\mathsf{S}:=\{1,\dots,S\}$, where ${N}_j$ is the number of samples in the $j^{\mathrm{th}}$ scale, such that the total number of cycles is $N_\mathrm{cyc}=\prod_{j=1}^SN_j$. The pseudo-time instants $t_{n}$ can then be recovered as\footnote{Empty sums and empty products, e.g., $\sum_{i=a}^b\Box_i$ and $\prod_{i=a}^b\Box_i$ for $b<a$, take values of $0$ and $1$, respectively.}
\begin{equation}
t_{n}=T\sum\limits_{j\in\mathsf{S}}({n}_j-1)\prod\limits_{i=1}^{j-1}{N}_i+\tau_h, \quad \text{with} \quad n = (N_\tau-1)\sum\limits_{j\in\mathsf{S}}({n}_j-1)\prod\limits_{i=1}^{j-1}{N}_i+h.
\label{multitemp}
\end{equation}
Figure~\ref{fig:tsep} shows a schematic representation of the multi-temporal scheme, where, for illustrative purposes, the time evolution of a generic scalar-valued variable $\chi$ is considered with 3 large-time scales. In this setting, we are able to distinguish processes occurring at multiple time scales, at the cost of introducing the \emph{course of dimensionality}. In particular, a quantity $\chi$ changes in a 3D spatial setting from $$\underbrace{\vphantom{\sum}\color{black}(\bm{x},n)\mapsto\chi_n(\bm{x})}_{\color{black}{{4} \text{ dimensions}}} \quad \text{to} \quad \underbrace{\vphantom{\sum}\color{black}(\bm{x},h,\mathfrak{n})\mapsto\chi_{h\mathfrak{n}}(\bm{x})}_{\color{black}{{4+S} \text{ dimensions}}}, \quad \text{with}  \quad \mathfrak{n}:=\{n_1,\dots,n_S\}.$$

The PGD framework~\cite{chinesta2014} approaches this problem by seeking solutions admitting a low-rank approximation. In the present study, we consider the following ansatz for the displacement field:
\begin{equation}
\bm{u}_{n}(\bm{x})\approx\bm{u}^{\mathrm{pgd}}_{h\mathfrak{n}}(\bm{x}) \coloneqq \sum_{i=1}^{M}{\color{black}\bm{{\varphi}}^{(i)}_h(\bm{x})}\prod_{j=1}^S{\color{black}\vartheta^{(j)(i)}_{n_j}}.
\label{eq:PGD1}
\end{equation}
For each mode $i$, $\bm{{\varphi}}^{(i)}\colon\Omega\times\mathsf{T}\to\mathbb{R}^d$ is an intra-cycle displacement function that evolves in the small-time interval $[0,T]$, while $\vartheta^{(j)(i)}\colon\mathsf{N}_j\to\mathbb{R}$ is a discrete function of $n_j$ that only varies in the large-time scale $j\in\mathsf{S}$, modulating the intra-cycle displacements. The ansatz~\eqref{eq:PGD1} is assumed to converge to the solution of the mechanical balance equation~\eqref{eq:st_adm1} as $M\to
\infty$. However, in practice, a finite number of modes are expected to render a reasonable approximation. Note that in contrast with a posteriori MOR, the modes are a priori unknown. Then, assuming all modes from $i=1$ to $i=m-1$ are known, the goal is to find the unknown enrichment functions $\bm{\varphi}^{(m)}_h(\bm{x})$ and $\vartheta^{(j)(m)}_{{n_j}}$ for all $\bm{x}\in\Omega$, $h\in\mathsf{T}$, and  $n_j\in\mathsf{N}_j\,\forall\,j\in\mathsf{S}$. 

It is worth noting that the space-time variation of $\bm{\varphi}^{(m)}_h(\bm{x})$, describing a characteristic cycle and possessing only $N_\tau$ temporal degrees of freedom, is expected to be less complex than that of the full solution $\bm{u}^{\mathrm{pgd}}_{h\mathfrak{n}}(\bm{x})$. Moreover, in the case of a single large-time scale, the one-dimensional function $\vartheta^{(1)(m)}_{{n_1}}$ possesses as many degrees of freedom as the number of cycles $N_1=N_\mathrm{cyc}$. By separating this function into multiple large-time scales, we may seek simpler functions $\vartheta^{(j)(m)}_{{n_j}}$ with $N_j$ degrees of freedom, such that the total number of large-time degrees of freedom changes from $N_\mathrm{cyc}=\prod_{j=1}^SN_j$ to $\sum_{j=1}^SN_j$.

Let us first define, in view of equation~\eqref{eq:constitutive}, the stress approximation
\begin{equation}
\sig_{h\mathfrak{n}}(\bm{x})\approx\frac{\partial\psi\big(\nabla^\mathrm{s}\bm{u}^{\mathrm{pgd}}_{h\mathfrak{n}},\mathbf{a}_{h\mathfrak{n}}\big)}{\partial\nabla^\mathrm{s}\bm{u}^{\mathrm{pgd}}_{h\mathfrak{n}}}(\bm{x}), \quad \nabla^\mathrm{s}\bm{u}^{\mathrm{pgd}}_{h\mathfrak{n}}(\bm{x}) = \sum_{i=1}^{M}{\color{black}\nabla^\mathrm{s}\bm{{\varphi}}^{(i)}_h(\bm{x})}\prod_{j=1}^S{\color{black}\vartheta^{(j)(i)}_{n_j}}.
\label{eq:stressansatz}
\end{equation}
Next, letting , $\tilde{\bm{\varphi}}_h\in\mathscr{U}_0$ and $\tilde{{\vartheta}}_{n_j}^{(j)}\in\mathbb{R}$, we define a virtual displacement field in separated form:
\begin{equation}
\begin{aligned}
\tilde{\bm{u}}^{\mathrm{pgd}}_{h\mathfrak{n}}(\bm{x}):=\tilde{\bm{\varphi}}_h(\bm{x})\prod_{j\in\mathsf{S}}{\vartheta}^{(j)(m)}_{n_j} + \bm{\varphi}^{(m)}_h(\bm{x})\sum_{j=1}^S\tilde{\vartheta}_{n_j}^{(j)}\prod_{k\in\mathsf{S}\setminus\{j\}}{\vartheta}^{(k)(m)}_{n_k}.
\label{eq:PGD_virt}
\end{aligned}
\end{equation}
Using expressions~\eqref{eq:PGD1}--\eqref{eq:PGD_virt}, and integrating in $[0,T]$ numerically, equation~\eqref{eq:st_adm1} takes the form
\begin{equation}
\sum_{h=1}^{N_\tau}\sum_{n_1=1}^{N_1}\dots \sum_{n_S=1}^{N_S}w_{h}\bigg[\int_{\Omega}\sig_{h\mathfrak{n}}:\nabla^\mathrm{s}\tilde{\bm{u}}^{\mathrm{pgd}}_{h\mathfrak{n}}\,\mathrm{d}\bm{x}-\int_{\Omega}\boldsymbol{b}_{h\mathfrak{n}}\cdot\tilde{\bm{u}}^{\mathrm{pgd}}_{h\mathfrak{n}}\,\mathrm{d}\bm{x} - \int_{\Gamma_\mathrm{N}}\bar{\bm{t}}_{h\mathfrak{n}}\cdot\tilde{\bm{u}}^{\mathrm{pgd}}_{h\mathfrak{n}}\,\mathrm{d}S\bigg]=0,
\label{eq:pgd1}
\end{equation}
where $w_h$ is a weight factor. This expression yields the following coupled problem: find $\bm{\varphi}_h^{(m)}\in\mathscr{U}_{{0}}$ and $\vartheta^{(j)(m)}_{n_j}\in\mathbb{R}$ for all $h\in\mathsf{T}$, $n_j\in\mathsf{N}_j$, and $j\in\mathsf{S}$ such that:
\begin{equation}
\sum_{h=1}^{N_\tau}\sum_{n_1=1}^{N_1}\dots \sum_{n_S=1}^{N_S}w_{h}\bigg[\int_{\Omega}\big(\sig_{h\mathfrak{n}}:\nabla^\mathrm{s}\tilde{\bm{\varphi}}_h-\bm{b}_{h\mathfrak{n}}\cdot\tilde{\bm{\varphi}}_h\big)\,\mathrm{d}\bm{x} -  \int_{\Gamma_\mathrm{N}}\bar{\bm{t}}_{h\mathfrak{n}}\cdot\tilde{\bm{\varphi}}_h\,\mathrm{d}S\bigg]\prod_{j\in\mathsf{S}}{\vartheta}^{(j)(m)}_{n_j}=0,
\label{eq:pgd_spat}
\end{equation}
\begin{equation}
\sum_{h=1}^{N_\tau}\sum_{n_1=1}^{N_1}\dots \sum_{n_S=1}^{N_S}w_{h}\bigg[\int_{\Omega}\big(\sig_{h\mathfrak{n}}:\nabla^{\mathrm{s}}\bm{\varphi}_h^{ (m)}-\bm{b}_{h\mathfrak{n}}\cdot\bm{\varphi}^{(m)}_h\big)\,\mathrm{d}\bm{x} - \int_{\Gamma_\mathrm{N}}\boldsymbol{t}_{h\mathfrak{n}}\cdot\bm{\varphi}^{(m)}_h\mathrm{d}S\bigg]\sum_{j=1}^S\tilde{\vartheta}_{n_j}^{(j)}\prod_{k\in\mathsf{S}\setminus\{j\}}{\vartheta}^{(k)(m)}_{n_k}=0.
\label{eq:pgd_temp}
\end{equation}
This system is non-linear due to the multiplicative decomposition~\eqref{eq:PGD1}, even in the case of linear elasticity. Its solution thus requires an iterative process, typically consisting of a fixed-point iteration scheme. 

If we consider, for now, that the time history of internal variables $\{\mathbf{a}_n\}_{n\in\mathsf{I}}$ is fixed, the system~\eqref{eq:pgd_spat}--\eqref{eq:pgd_temp}  allows for a reduced number of PDE solutions. In particular, the PDE~\eqref{eq:st_adm1} need not be solved at each $n\in\mathsf{I}$. Instead, fixing $\{\vartheta^{(j)(m)}_{n_j}\}_{n_j\in\mathsf{N}_j}$ for all $j\in\mathsf{S}$, the PDE~\eqref{eq:pgd_spat} is solved for $\bm{\varphi}_h^{(m)}$ at each $h\in\mathsf{T}$. Then, fixing $\{\bm{\varphi}_h^{(m)}\}_{h\in\mathsf{T}}$, equation~\eqref{eq:pgd_temp} renders a system of $S$  algebraic equations, one for each large-time scale $j\in\mathsf{S}$, to be solved for $\vartheta^{(j)(m)}_{n_j}$ at each $n_j\in\mathsf{N}_j$. If this procedure is iteratively performed for each mode until convergence, the total number of PDE solutions changes from $N_t= (N_\tau-1) N_\mathrm{cyc}$ in the classical incremental scheme to $N_\tau \sum_{i=1}^{M}N_\mathrm{iter}^{(i)}$ in the PGD formulation, where $N_\mathrm{iter}^{(i)}$ is the number of iterations required to solve~\eqref{eq:pgd_spat}--\eqref{eq:pgd_temp} for mode $i$. Hence, a lower number of PDE solutions can be achieved for large cycle numbers, provided that limited number of modes yields a good approximation, with a limited number of iterations per mode.

It bears emphasis that these observations assume a given time history of internal variables $\{\mathbf{a}_n\}_{n\in\mathsf{I}}$. However, $\{\mathbf{a}_n\}_{n\in\mathsf{I}}$ follows from the solution of the evolution equation~\eqref{eq:biot} in $\Omega\times\mathrm{I}$, which is coupled to~\eqref{eq:pgd_spat}--\eqref{eq:pgd_temp} through the displacement field. As a consequence, equations~\eqref{eq:pgd_spat} and~\eqref{eq:pgd_temp} are non-linear individually in $\bm{\varphi}_h^{(m)}$ and $\vartheta^{(j)(m)}_{n_j}$. Moreover, due to the incremental nature of~\eqref{eq:biot}, the internal variables are not separable into multiple time scales, such that $\sig_{h\mathfrak{n}}$,  related to $\nabla^\mathrm{s}\bm{u}^{\mathrm{pgd}}_{h\mathfrak{n}}$ and $\mathbf{a}_{h\mathfrak{n}}$ in incremental form, is also not separable. The numerical approach to this problem is discussed in section~\ref{sec:numimp}. 

\begin{remark}
    The present work focuses on the separated solution of the mechanical balance system~\eqref{eq:pgd_spat}--\eqref{eq:pgd_temp} along with a conventional solution of the evolution equation~\eqref{eq:biot}. Further developments, focusing on specific models with internal variables, are required to establish a separated solution of the evolution equation, which would allow us to fully exploit the computational advantages of the PGD formalism. We leave these developments to future studies. Further details on this aspect are given in section~\ref{sec:numimp}. 
\end{remark}

\begin{remark}
The present formulation does not separate the small-time coordinates $\tau_h$ from the spatial coordinates $\bm{x}$. This specific form was chosen because the internal variables $\mathbf{a}$ are expected to evolve rapidly within each cycle. A straightforward modification may be considered to render a fully separated form, where each mode in~\eqref{eq:PGD1} is rewritten as $\bm{\varphi}^{(i)}(\bm{x})\gamma^{(i)}_h\prod_{j=1}^S{\vartheta^{(j)(i)}_{n_j}}$, with $\gamma\colon\mathsf{T}\to\mathbb{R}$ a scalar function that governs the small-time scale. An additional scalar equation, to be solved for this function, is then appended to the coupled system~\eqref{eq:pgd_spat}--\eqref{eq:pgd_temp}. The advantage is that $\bm{\varphi}^{(i)}$ becomes time-independent, with a single PDE solution required per iteration, per mode.  This option may be considered in future~works.
\end{remark}

\section{Application to ratcheting plasticity and numerical implementation}

This section presents an application of the formulation in section~\ref{sec:form} to cyclic plasticity. Specifically, we focus on the numerical implementation of the multi-temporal system~\eqref{eq:pgd_spat}--\eqref{eq:pgd_temp} coupled with a ratcheting plasticity model. This model was originally proposed by~\citet{houlsby2017} for monopile foundations and further extended in~\citet{ulloa2021b} to encompass general elastoplastic solids and ductile failure. It is particularly interesting for the present work, as ratcheting entails a non-periodic response to periodic loading, exhibiting, for instance, cyclic strain accumulations under constant-amplitude~stress. 

\subsection{Constitutive model}
\label{sec:const}

Let us briefly present the constitutive model according to the formulation outlined in section~\ref{sec:setup} and~\citet{ulloa2021a}. We consider the set of internal variables $\mathbf{a}\coloneqq \{\bm{\varepsilon}^\mathrm{p},\kappa,\bm{\varepsilon}^\mathrm{r}\}$, where $\epsp\colon\Omega\times\setI\to\mathbb{R}^{3\times 3}_\mathrm{sym}$ is the plastic strain tensor, $\kappa\colon\Omega\times\setI\to\mathbb{R}_+$ is an isotropic hardening variable, and $\bm{\varepsilon}^\mathrm{r}\colon\Omega\times\setI\to\mathbb{R}^{3\times 3}_\mathrm{sym}$ is a ratcheting strain tensor. The corresponding set of generalized stresses~\eqref{eq:constitutive}$_2$ is given by $\mathbf{s}\coloneqq\{\bm{s}^\mathrm{p},h,\bm{s}^\mathrm{r}\}$. 

The free energy density is defined within the small strain setting as 
\begin{equation}
\psi(\bm{\varepsilon},\mathbf{a})\coloneqq\frac{1}{2}\big(\bm{\varepsilon}-\bm{\varepsilon}^\mathrm{p}-\bm{\varepsilon}^\mathrm{r}\big):\bm{\mathsf{C}}:\big(\bm{\varepsilon}-\bm{\varepsilon}^\mathrm{p}-\bm{\varepsilon}^\mathrm{r}\big)+\frac{1}{2}H^\mathrm{kin}\epsp:\epsp+\frac{1}{2}H^\mathrm{iso}\kappa^2,
\end{equation}
where $\bm{\mathsf{C}}$ is the elasticity tensor,  $H^\mathrm{kin}$ is the kinematic hardening modulus, and $H^\mathrm{iso}$ is the isotropic hardening modulus. On the other hand, the dissipation potential reads

\begin{equation}
    \phi(\dot{\mathbf{a}};\mathbf{s})\coloneqq\begin{dcases}
\bigg(\sqrt{\frac{2}{3}}\sigma^\mathrm{p}+{\beta\Vert\bm{s}^\mathrm{r}_\mathrm{dev}\Vert}\bigg)\Vert \dot{\bm{\varepsilon}}^\mathrm{p}_\mathrm{dev}\Vert \quad \text{if} \  \ \begin{dcases} \dot{\kappa}\geq \sqrt{\frac{2}{3}}\Vert \dot{\bm{\varepsilon}}^\mathrm{p}_\mathrm{dev}\Vert, \   \Vert \dot{\bm{\varepsilon}}^\mathrm{r}_\mathrm{dev}\Vert\leq\beta\Vert \dot{\bm{\varepsilon}}^\mathrm{p}_\mathrm{dev}\Vert, \\  \tr \dot{\bm{\varepsilon}}^\mathrm{p}=0, \ \tr \dot{\bm{\varepsilon}}^\mathrm{r}=0, \end{dcases}\\
+\infty \quad \text{otherwise},
\end{dcases}
\end{equation}
where $\beta\in[0,1]$ is a ratcheting constant. Using equations~\eqref{eq:constitutive} and~\eqref{eq:biot}, it can be shown~\cite{ulloa2021a} that these energy functions yield a constitutive model with the following ingredients:
\begin{align}
&\text{Stress-strain relation} &&\sig=\bm{\mathsf{C}}:\big(\bm{\varepsilon}-\bm{\varepsilon}^\mathrm{p}-\bm{\varepsilon}^\mathrm{r}\big);\label{eq:constitutitve2}   \\
&\text{Generalized stresses}
&&\bm{s}^\mathrm{p}=\sig - H^\mathrm{kin}\bm{\varepsilon}^\mathrm{p}, \quad h=- H^{\mathrm{iso}}\kappa,\quad \bm{s}^\mathrm{r}=\sig;  \\ 
&\text{Yield function} &&f(\bm{s}^\mathrm{p},h)=\Vert\bm{s}^\mathrm{p}_{\mathrm{dev}}\Vert-\sqrt{\frac{2}{3}}(\sigma^\mathrm{p}+h); \\
&\text{Plastic potential} &&g(\bm{s}^\mathrm{p},h,\bm{s}^\mathrm{r})=f(\bm{s}^\mathrm{p},h)+\beta\Vert\bm{s}^\mathrm{r}_\mathrm{dev}\Vert; \\
&\text{KKT conditions} &&f\leq 0,\quad \dot{\lambda} f = 0, \quad \dot{\lambda}\geq 0;\\
&\text{Flow rule} && (\dot{\bm{\varepsilon}}^{\mathrm{p}},\dot{\kappa},\dot{\bm{\varepsilon}}^{\mathrm{r}})\in\dot{\lambda}\,\partial g(\bm{s}^\mathrm{p},h,\bm{s}^\mathrm{r}).\label{eq:flowrule}
\end{align}
These equations are defined locally in $\Omega\times\mathrm{I}$ and are coupled to the system~\eqref{eq:pgd_spat}--\eqref{eq:pgd_temp}.
	
\subsection{Numerical implementation}
\label{sec:numimp}

\paragraph{Finite-element discretization} 
 Using standard finite-element notation, we perform spatial interpolation as $\bm{\varphi}_h(\bm{x})\approx\mathbf{N}(\bm{x})\discrete{\bm{\varphi}}_h$ and $\{\nabla^{\mathrm{s}}\bm{\varphi}_{h}\}(\bm{x})\approx \mathbf{B}(\bm{x})\discrete{\bm{\varphi}}_h$, with $\discrete{\bm{\varphi}}_h\in\mathbb{R}^{N_\mathrm{d}}$, where $N_\mathrm{d}$ is the number of spatial degrees of freedom. Hereinafter, $\{\Box\}$ is the Voigt representation of a second-order tensor $\Box$. Defining the external force vector as
\begin{equation}
\mathbf{f}^\mathrm{ext}_{h\mathfrak{n}}:= \int_{\Omega}\mathbf{N}^\mathrm{T}\bm{b}_{h\mathfrak{n}}\,\mathrm{d}\bm{x} + \int_{\Gamma_\mathrm{N}}\mathbf{N}^\mathrm{T}\bar{\bm{t}}_{h\mathfrak{n}}\,\mathrm{d}S, 
\end{equation}
the separated system~\eqref{eq:pgd_spat}--\eqref{eq:pgd_temp} takes the discrete residual form
\begin{equation}
\bm{R}^\varphi\big[\discrete{\bm{\varphi}}^{(m)}_h\big]\coloneqq\sum_{n_1=1}^{N_1}\dots \sum_{n_S=1}^{N_S}\prod_{j\in\mathsf{S}}{\vartheta}^{(j)(m)}_{n_j}\bigg[\int_\Omega\mathbf{B}^\mathrm{T}\{\sig_{h\mathfrak{n}}\}\,\mathrm{d}\bm{x}-\mathbf{f}^{\mathrm{ext}}_{h\mathfrak{n}}\bigg]=\bm{0} \quad \forall\,h\in\mathsf{T},
\label{eq:pgd_spat_disc}
\end{equation}
\begin{equation}
\begin{aligned}
{R}^{\vartheta}\big[\vartheta^{(j)(m)}_{n_j}\big]\coloneqq\sum_{h=1}^{N_\tau}\sum_{\substack{{n_k=1} \\ \forall k\in\mathsf{S}\setminus\{j\}}}^{N_k}\prod_{\ell\in\mathsf{S}\setminus\{j\}}\hspace{-0.5em}{\vartheta}^{(\ell)(m)}_{n_\ell}w_{h}\,\discrete{\bm{\varphi}}_h^{(m)\mathrm{T}}\bigg[\int_\Omega\mathbf{B}^\mathrm{T}\{\sig_{h\mathfrak{n}}\}\,\mathrm{d}\bm{x}-\mathbf{f}^{\mathrm{ext}}_{h\mathfrak{n}}\bigg]=0  \quad \forall\,n_j\in\mathsf{N}_j.
\end{aligned}
\label{eq:pgd_temp_disc}
\end{equation}
In view of equation~\eqref{eq:constitutitve2}, the Voigt form of the stress tensor reads 
\begin{equation}
    \{\sig_{h\mathfrak{n}}\}=\mathbf{D}\,\big\{\nabla^ \mathrm{s}\bm{u}^{\mathrm{pgd}}_{h\mathfrak{n}}-\epsp_{n}-\epsr_{n}\big\}, 
    \label{eq:constitutive3}
\end{equation}
where $\mathbf{D}$ is the matrix representation of the elasticity tensor $\bm{\mathsf{C}}$. The solution of the discrete system~\eqref{eq:pgd_spat_disc}--\eqref{eq:pgd_temp_disc} allows us to reconstruct the full-time history of nodal displacements $\{\discrete{\bm{\mathrm{u}}}_n\}_{n\in\mathsf{I}}\in\mathbb{R}^{ N_\mathrm{d}\times N_t}$ as
\begin{equation}
\{\discrete{\bm{\mathrm{u}}}_n\}_{n\in\mathsf{I}}=\sum_{i=1}^M\big\{{\discrete{\bm{\mathrm{\varphi}}}}_h^{(i)}\big\}_{h\in\mathsf{T}}\otimes\big\{\vartheta_{n1}^{(1)(i)}\big\}_{n_1\in\mathsf{N}_1}\otimes\dots\otimes{{\big\{\vartheta}}_{nS}^{(S)(i)}\big\}_{n_S\in\mathsf{N}_S}.
\label{eq:PGD2}
\end{equation}
Consequently, the total number of displacement degrees of freedom in space-time changes from $N_\mathrm{d}N_t=N_\mathrm{d}(N_\tau-1) N_\mathrm{cyc}= N_\mathrm{d}(N_\tau-1)\prod_{j=1}^SN_j$ in the classical incremental scheme to $MN_\mathrm{d}N_\tau+M\sum_{j=1}^SN_j$ in the PGD formulation, representing a possibly major reduction for a finite mode number $M$. 

\paragraph{Fixed-point iterations}
The system~\eqref{eq:pgd_spat_disc}--\eqref{eq:pgd_temp_disc} is solved for each modal enrichment in a fixed-point iteration scheme. Specifically, with all modes known up to $m-1$,  equation~\eqref{eq:pgd_spat_disc} is solved incrementally for $\{\discrete{\bm{\varphi}}_h^{(m)}\}_{h\in\mathsf{T}}$ with fixed large-time function values $\{\vartheta^{(j)(m)}_{n_j}\}_{n_j\in\mathsf{N}_j} \ \forall\,j\in\mathsf{S}$. Then, fixing the intra-cycle displacements $\{\discrete{\bm{\varphi}}_h^{(m)}\}_{h\in\mathsf{T}}$, the $S$ algebraic equations in~\eqref{eq:pgd_temp_disc} are solved for each $\{\vartheta^{(j)(m)}_{n_j}\}_{n_j\in\mathsf{N}_j}$, and the procedure is repeated iteratively until convergence, that is, until the $L^2$ norm of the algorithmic corrections is sufficiently small for the intra-cycle displacements and all the large-time function values.

\paragraph{Linearization}

The main difficulty in this procedure lies in the fact that the stress tensor~\eqref{eq:constitutive3} is not a priori separable into multiple time scales, with $\epsp_{n}$ and $\epsr_{n}$ only known from the incremental solution of the constitutive model~\eqref{eq:constitutitve2}--\eqref{eq:flowrule}. Let us consider two possible approaches.

\begin{enumerate}[wide, labelwidth=!, labelindent=0pt]
    \item \emph{Consistent linearization.}   Equations~\eqref{eq:pgd_spat_disc} and~\eqref{eq:pgd_temp_disc} may be solved individually with a consistent Newton-Raphson scheme. In this way, with fixed  $\{\vartheta^{(j)(m)}_{n_j}\}_{n_j\in\mathsf{N}_j} \ \forall\,j\in\mathsf{S}$,~\eqref{eq:pgd_spat_disc} is linearized as
    \begin{equation}
    -\sum_{n_1=1}^{N_1}\dots \sum_{n_S=1}^{N_S}\prod_{j\in\mathsf{S}}\left({\vartheta}^{(j)(m)}_{n_j}\right)^2\hspace{-0.25em}\,\mathbf{K}^{\mathrm{tan}(\alpha)}_{h\mathfrak{n}}\,\Delta\discrete{\bm{\varphi}}^{(m)(\alpha+1)}_h=\bm{R}^\varphi\big[\discrete{\bm{\varphi}}^{(m)(\alpha)}\big] \quad \forall\,h\in\mathsf{T},
        \label{eq:pgd_spat_disc_lin}
    \end{equation}
    where the superscript $(\alpha)$ is an iteration counter and $\Delta\Box^{(\alpha+1)}\coloneqq\Box^{(\alpha+1)}-\Box^{(\alpha)}$ is an algorithmic increment. The tangent stiffness matrix $\mathbf{K}^\mathrm{tan}_{h\mathfrak{n}} \in \mathbb{R}^{N_\mathrm{d}\times N_\mathrm{d}}$ takes the form
    \begin{equation}
    \mathbf{K}^\mathrm{tan}_{h\mathfrak{n}}=\int_{\Omega}{\mathbf{B}(\bm{x})}^{\mathrm{T}}\mathbf{D}^\mathrm{tan}_{h\mathfrak{n}}(\bm{x})\,{\mathbf{B}(\bm{x})}\,\mathrm{d}\bm{x},
    \label{eq:Ktan}
    \end{equation}
    where, in view of equation~\eqref{eq:constitutive3} and the flow rule~\eqref{eq:flowrule}, the consistent tangent operator reads
    \begin{equation} \mathbf{D}^\mathrm{tan}_{h\mathfrak{n}}=\frac{\partial \{\sig_{h\mathfrak{n}}\}}{\partial\{\nabla^ \mathrm{s}\bm{u}^{\mathrm{pgd}}_{h\mathfrak{n}}\}}\equiv\mathbf{D}\bigg[\mathbf{I}-\frac{\partial}{\partial\{\eps_n\}}\Delta\lambda_{n}\{\hat{\bm{n}}_{n}+\beta\hat{\bm{r}}_{n}\}\bigg], \quad \Delta\lambda_{n}\coloneqq\lambda_{n}-\lambda_{n-1},
    \end{equation}
    with $\hat{\bm{n}}\coloneqq\partial\Vert\bm{s}^\mathrm{p}_{\mathrm{dev}}\Vert$ and $\hat{\bm{r}}\coloneqq\partial{\Vert\bm{s}^\mathrm{r}_{\mathrm{dev}}\Vert}$. Equation~\eqref{eq:pgd_spat_disc_lin} is solved iteratively for all $\discrete{\bm{\varphi}}_h^{(m)}$ until convergence, with $\{\epsp_{n}\}_{n\in\mathsf{I}}$ and $\{\epsr_{n}\}_{n\in\mathsf{I}}$ updated at each iteration from the incremental solution of the constitutive model~\eqref{eq:constitutitve2}--\eqref{eq:flowrule}. On the other hand, with fixed $\{\discrete{\bm{\varphi}}_h^{(m)}\}_{h\in\mathsf{T}}$, equation~\eqref{eq:pgd_temp_disc} is linearized as
    \begin{equation}
        -\sum_{h=1}^{N_\tau}\sum_{\substack{{n_k=1} \\ \forall k\in\mathsf{S}\setminus\{j\}}}^{N_k}\hspace{-0.25em}\prod_{\ell\in\mathsf{S}\setminus\{j\}}\hspace{-0.5em}\left({\vartheta}^{(\ell)(m)}_{n_\ell}\right)^2\hspace{-0.5em}w_h\,{\discrete{\bm{\varphi}}^{(m)\mathrm{T}}_h}\mathbf{K}^{\mathrm{tan}(\beta)}_{h\mathfrak{n}}\,\discrete{\bm{\varphi}}_h^{(m)}\,\Delta\vartheta^{(j)(m)(\beta+1)}_{n_j}={R}^{\vartheta}\big[\vartheta^{(j)(m)(\beta)}_{n_j}\big]\quad \forall\,n_j\in\mathsf{N}_j,
        \label{eq:pgd_temp_disc_lin}
    \end{equation}
    where the superscript $(\beta)$ is the iteration counter. This equation is solved iteratively for all $\vartheta^{(j)(m)}_{n_j}$ until convergence, with $\{\epsp_{n}\}_{n\in\mathsf{I}}$ and $\{\epsr_{n}\}_{n\in\mathsf{I}}$ also updated at each iteration from~\eqref{eq:constitutitve2}--\eqref{eq:flowrule}.
    
    \item \emph{Decoupled scheme.} The procedure above requires evaluating the full-time history of internal variables $\{\epsp_{n}\}_{n\in\mathsf{I}}$ and $\{\epsr_{n}\}_{n\in\mathsf{I}}$ and the consistent tangent~\eqref{eq:Ktan} for each mode, at each fixed-point iteration of the system~\eqref{eq:pgd_spat_disc}--\eqref{eq:pgd_temp_disc}, at each Newton-Raphson iteration of~\eqref{eq:pgd_spat_disc_lin}, and at each Newton-Raphson iteration of~\eqref{eq:pgd_temp_disc_lin}. A way to alleviate these heavy computations is to consider a single Newton-Raphson iteration for both~\eqref{eq:pgd_spat_disc_lin} and~\eqref{eq:pgd_temp_disc_lin}. This approach amounts to resolving the full-time history of the displacement field by solving the fixed-point iterations of~\eqref{eq:pgd_spat_disc}--\eqref{eq:pgd_temp_disc}  with fixed time histories $\{\epsp_{n}\}_{n\in\mathsf{I}}$ and $\{\epsr_{n}\}_{n\in\mathsf{I}}$. As such, these equations are now individually linear in $\discrete{\bm{\varphi}}^{(m)}_h$ and ${\vartheta}^{(j)(m)}_{n_j}$ and take the form
    \begin{equation}
    \begin{aligned}
    \sum_{n_1=1}^{N_1}\dots \sum_{n_S=1}^{N_S}\prod_{j\in\mathsf{S}}{\vartheta}^{(j)(m)}_{n_j}\Bigg(\prod_{j\in\mathsf{S}}{\vartheta}^{(j)(m)}_{n_j}\mathbf{K}\,\discrete{\bm{\varphi}}^{(m)}_h&+\mathbf{K}\sum_{i=1}^{m-1}\discrete{\bm{\varphi}}^{(i)}_h\prod_{j\in\mathsf{S}}{\vartheta}^{(j)(i)}_{n_j}\\&\qquad-\int_\Omega\mathbf{B}^\mathrm{T}\mathbf{D}\,\{\epsp_{h\mathfrak{n}}+\epsr_{h\mathfrak{n}}\}\,\mathrm{d}\bm{x}-\mathbf{f}^{\mathrm{ext}}_{h\mathfrak{n}}\Bigg)=\bm{0} \quad \forall\,h\in\mathsf{T},
    \end{aligned}
    \label{eq:pgd_spat_disc_lin0}
    \end{equation}
    \begin{equation}
    \begin{aligned}
    \sum_{h=1}^{N_\tau}\sum_{\substack{{n_k=1} \\ \forall k\in\mathsf{S}\setminus\{j\}}}^{N_k}\hspace{-0.25em}\prod_{\ell\in\mathsf{S}\setminus\{j\}}\hspace{-0.5em}{\vartheta}^{(\ell)(m)}_{n_\ell}w_{h}\,{\discrete{\bm{\varphi}}_h^{(m)\mathrm{T}}}\Bigg(\mathbf{K}\,\discrete{\bm{\varphi}}^{(m)}_h\hspace{-0.75em}\prod_{\ell\in\mathsf{S}\setminus\{j\}}\hspace{-0.5em}&{\vartheta}^{(\ell)(m)}_{n_\ell}\,{\vartheta}^{(j)(m)}_{n_j} + \mathbf{K}\sum_{i=1}^{m-1}\discrete{\bm{\varphi}}^{(i)}_h\prod_{j\in\mathsf{S}}{\vartheta}^{(j)(i)}_{n_j}  \\
    &  -\int_\Omega\mathbf{B}^\mathrm{T}\mathbf{D}\,\{\epsp_{h\mathfrak{n}}+\epsr_{h\mathfrak{n}}\}\,\mathrm{d}\bm{x}-\mathbf{f}^{\mathrm{ext}}_{h\mathfrak{n}}\Bigg)=0\quad \forall\,n_j\in\mathsf{N}_j \hspace*{\fill},
    \end{aligned}
    \label{eq:pgd_temp_disc_lin0}
    \end{equation}
    with the constant stiffness matrix $\mathbf{K} \in \mathbb{R}^{N_\mathrm{d}\times N_\mathrm{d}}$ given by 
    \begin{equation}
\mathbf{K}=\int_{\Omega}{\mathbf{B}(\bm{x})}^{\mathrm{T}}\mathbf{D}(\bm{x})\,{\mathbf{B}(\bm{x})}\,\mathrm{d}\bm{x}.
    \label{eq:K}
    \end{equation}
    The solution of this system is then employed to update the full-time histories $\{\epsp_{n}\}_{n\in\mathsf{I}}$ and $\{\epsr_{n}\}_{n\in\mathsf{I}}$ from the constitutive model~\eqref{eq:constitutitve2}--\eqref{eq:flowrule}. The process is repeated iteratively until convergence.  
\end{enumerate}

We note that the consistent linearization procedure is rather cumbersome but may converge faster. On the other hand, the decoupled scheme is easy to implement but may result in slow convergence. \tB{Alternatively, a singular value decomposition of the time-dependent tangent stiffness matrix, as done in~\cite{arjoune2022} for geometrically non-linear systems with a single temporal scale, may be introduced in the present multi-temporal framework. Moreover, a partially consistent linearization may be considered, e.g., letting the tangent stiffness matrix evolve only in certain time scales. Nevertheless, for the sake of simplicity, we employ the decoupled scheme in the numerical simulations that follow.}

\paragraph{Normalization} It is convenient to normalize the modes as $\bm{u}^{\mathrm{pgd}}_{h\mathfrak{n}}(\bm{x})=\sum_{i=1}^{M}\zeta^{(i)}\hat{\bm{{\varphi}}}^{(i)}_h(\bm{x})\prod_{j=1}^S\hat{\vartheta}^{(j)(i)}_{n_j}$, with $\hat\Box$ denoting the $L^2$ normalization of a function $\Box$. Here, the magnitude of the coefficient $\zeta^{(i)}\in\mathbb{R}_+$ represents the contribution of mode $i$, with a small value indicating that the approximation has converged. The algorithmic procedure then consists of assuming $\zeta^{(m)}=1$ and solving for the modal enrichments, as discussed above, with the coefficients $\{\zeta^{(i)}\}_{i=1}^{m-1}$ known for the previous modes. Then, the solutions $\{\discrete{\bm{\varphi}}_h^{(m)}\}_{h\in\mathsf{T}}$ and $\{\vartheta^{(j)(m)}_{n_j}\}_{n_j\in\mathsf{N}_j} \, \forall\,j\in\mathsf{S}$ are normalized and the coefficients $\{\zeta^{(i)}\}_{i=1}^{m}$ are updated from the linear algebraic equation
\begin{equation}
     \sum_{h=1}^{N_\tau}\sum_{n_1=1}^{N_1}\dots \sum_{n_S=1}^{N_S}w_{h}\sum_{r=1}^{m}\prod_{k=1}^S\hat{\vartheta}^{(k)(r)}_{n_k}{\discrete{\hat{\bm{\varphi}}}^{(r)\mathrm{T}}_h}\Bigg(\mathbf{K}\sum_{i=1}^{m}\prod_{j=1}^S\discrete{\hat{\bm{\varphi}}}^{(i)}_h\hat{\vartheta}^{(j)(i)}_{n_j}\zeta^{(i)}-\int_\Omega\mathbf{B}^\mathrm{T}\mathbf{D}\,\{\epsp_{h\mathfrak{n}}+\epsr_{h\mathfrak{n}}\}\,\mathrm{d}\bm{x}-\mathbf{f}^{\mathrm{ext}}_{h\mathfrak{n}}\Bigg)=0.
\end{equation}
\begin{remark}
\label{rem:cont}
The formulation above does not enforce continuity in time, i.e., the displacements~\eqref{eq:PGD1} may exhibit temporal jumps from one cycle to another. However, as will be shown in the numerical examples, an approximately continuous solution can be achieved with a sufficient number of modes.
\end{remark}

\FloatBarrier
\section{Numerical simulations}
\label{sec:sims}

This section presents numerical simulations that demonstrate the feasibility of identifying multi-temporal solutions of the form~\eqref{eq:PGD2} in scenarios where periodic loading conditions yield a quasi-periodic response, governed by the ratcheting plasticity model of section~\ref{sec:const}. We aim at reproducing the reference incremental solution using the proposed multi-temporal PGD. To this end, we consider in both examples the decoupled numerical implementation scheme in section~\ref{sec:numimp}, solving the system~\eqref{eq:pgd_spat_disc_lin0}--\eqref{eq:pgd_temp_disc_lin0} via fixed point iterations and iterating the decomposed full-time solution with the constitutive model~\eqref{eq:constitutitve2}--\eqref{eq:flowrule}.

\subsection{Perforated plate}

\begin{figure}[h!]
\centering
\includegraphics[width=0.85\textwidth,trim={0cm 0cm 0cm 0cm},clip]{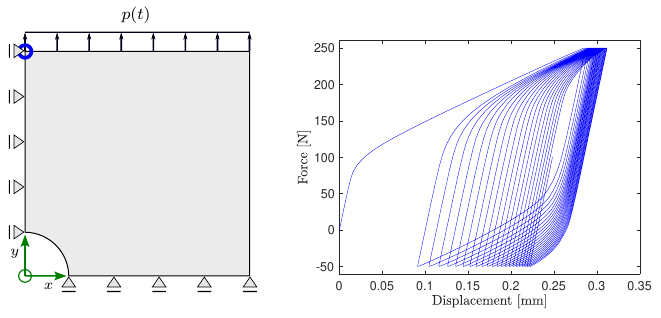}
\caption{Perforated plate under cyclic loading (left) and force-displacement response (right) at the top-left corner. The square specimen has a size of 30~mm by 30~mm while the circular notch has a radius of 6~mm. The plate is discretized with $N_\mathrm{d}=1640$ spatial degrees of freedom.}
\label{fig:plate1}
\end{figure}

The first example considers a perforated plate under cyclic loading and plane strain conditions (figure~\ref{fig:plate1}). The material obeys the ratcheting plasticity model of section~\ref{sec:const} with Young's modulus $E=205$~MPa, Poisson's ratio $\nu=0.3$, yield strength $\sigma^\mathrm{p}=100$~MPa, isotropic hardening modulus $H^\mathrm{iso}=1140$~MPa, kinematic hardening modulus $H^\mathrm{kin}=21640$~MPa, and ratcheting constant $\beta=0.4$. 

\begin{figure}[t!]
\centering
\includegraphics[width=1\textwidth,trim={0cm 0cm 0cm 0cm},clip]{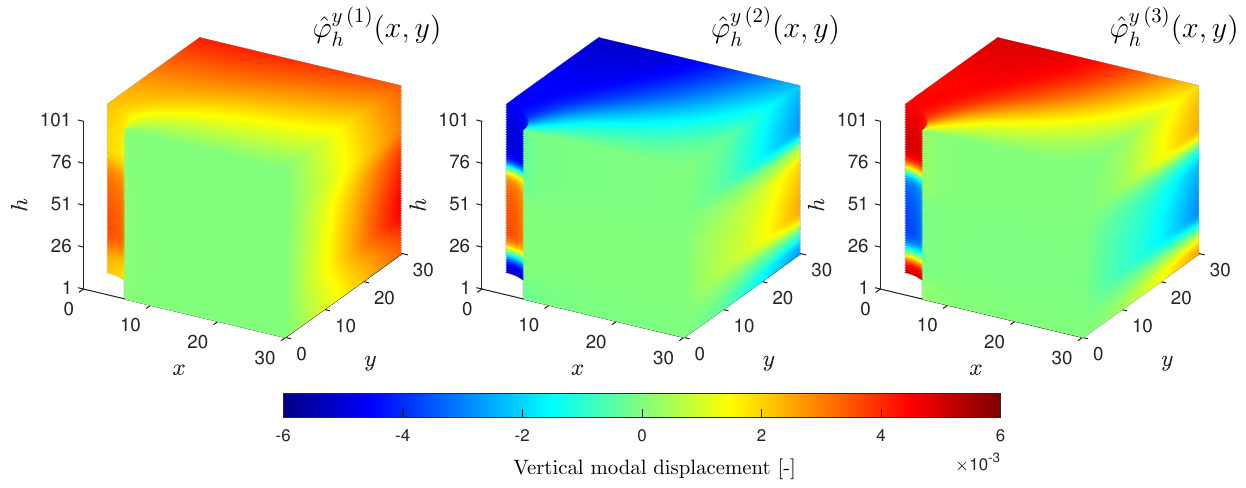}
\caption{Normalized intra-cycle displacements: spatial distribution in $(x,y)$ of the vertical ($y$) component and variation in the small-time scale $h$ for the first 3 modes.}
\label{fig:plate2}
\vspace{2em}
\includegraphics[width=1\textwidth,trim={0cm 0cm 0cm 0cm},clip]{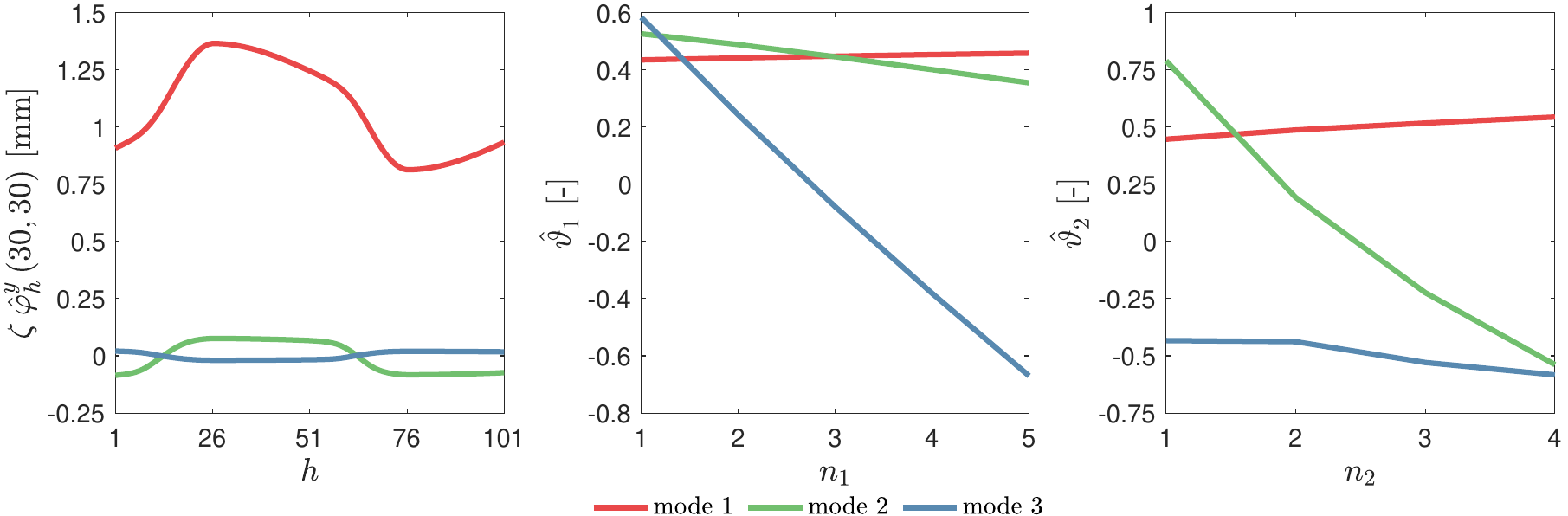}
\caption{Small-time variation at the top-left corner (left) and normalized first and second large-time function values (middle and right) for the first three modes.}
\label{fig:plate3}
\end{figure}

A constant-amplitude cyclic force is applied uniformly on the top side for 22 cycles, with the magnitude $p$ ranging from $-50$~N to $250$~N. Each cycle is discretized into 101 load steps. The force-displacement curve in figure~\ref{fig:plate2} shows a rather complex cyclic behavior, where ratcheting competes with cyclic hardening. In particular, the ratcheting response begins to stabilize after some cycles, with a pronounced  asymmetry in tension and compression.

\begin{figure}[h!]
\centering
\includegraphics[width=1\textwidth,trim={0cm 0cm 0cm 0cm},clip]{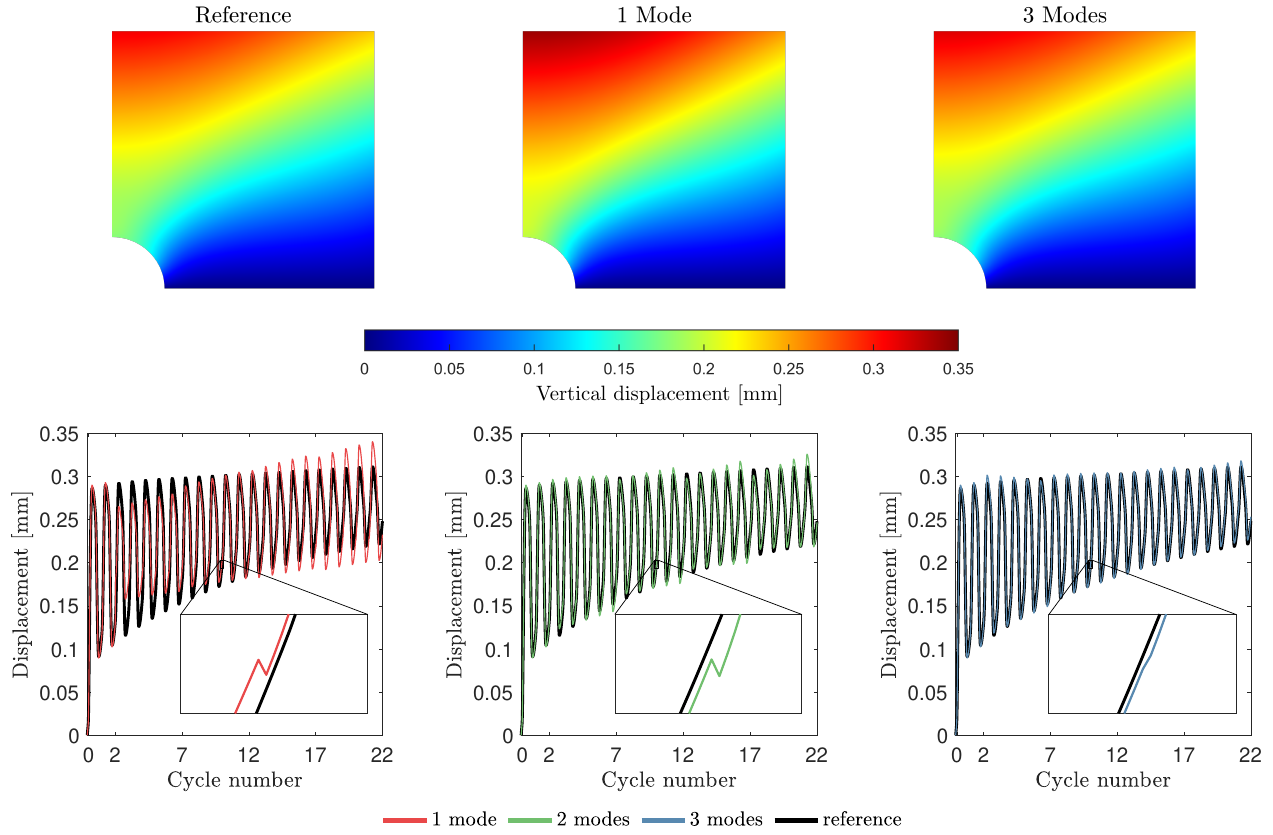}
\caption{Vertical displacements for the reference solution and the PGD solution with 1, 2, and 3 modes: spatial distribution (top) and time history at the top-left corner (bottom). PGD solution with 2 large-time scales ($N_1=5$ and $N_2=4)$.}
\label{fig:plate4}
\end{figure}

In this example, we resolve the first 2 cycles using the classical, single-scale incremental scheme. The response of the second cycle is then extrapolated periodically and employed as an initial guess to resolve the remaining 20 cycles using the proposed multi-temporal PGD. We consider 3 temporal scales: a small-time scale with $N_\tau=101$ and 2 large-time scales with $N_1=5$ and $N_2=4$. Although the number of simulated cycles is small, we already observe a considerable reduction of temporal degrees of freedom: assuming $M=3$, and noting that $N_\mathrm{d}=1640$, the total number of space-time degrees of freedom is reduced from $3.28\times10^6$ in the classical incremental scheme to $4.97\times10^5$ in the PGD~solution.

Figure~\ref{fig:plate2} shows the vertical component of the intra-cycle displacements $\{\discrete{\bm{\varphi}}_h^{(m)}\}_{h\in\mathsf{T}}$ for the first 3 modes. The first mode already shows the expected spatial variability, with the vertical displacements exhibiting a marked transition in an inclined region that starts at the notch, where the plastic strains concentrate. The second and third modes show a similar trend but change signs within the cycle, acting as correctors for the first mode. This result is better observed in figure~\ref{fig:plate3} (left), where the small-time variation scaled by the coefficients $\zeta^{(i)}$ is shown for a single point in space. The first mode clearly dominates the intra-cyclic response, with a modal coefficient $\zeta^{(1)}=222.68$, while the second and third modes have much lower magnitudes, with coefficients $\zeta^{(2)}=14.55$ and $\zeta^{(3)}=3.55$, respectively. The large-time functions in figure~\ref{fig:plate3} (middle and right) further show that most of the long-term variability, describing cyclic hardening, is captured by the second and third modes, while the first mode captures the competing ratcheting mechanism. The third mode seems to modulate the response in the first large-time scale, i.e., within a group of 5 cycles, while the second mode modulates the variability in the second large-time scale, i.e., between groups of 5~cycles.

Figure~\ref{fig:plate4} shows that the effect of the second and third modes in the long-term response is significant. The vertical displacements predicted by the PGD solution with 1, 2, and 3 modes are presented and compared with the reference incremental solution. While a single mode provides a reasonable spatial approximation, it exhibits some overestimation. Moreover, significant discrepancies are observed in the temporal evolution, transitioning from an underestimation in the initial cycles to an overestimation in the later cycles.  However, when using 2 modes, the response becomes much closer to the reference. Nevertheless, a close inspection reveals small temporal discontinuities from one cycle to another (see remark \ref{rem:cont}). Interestingly, after adding a third mode, the discontinuities are notably reduced and the overall reference response is approximated well. 

\FloatBarrier
\subsection{Monopile foundation}
 
\begin{figure}[t!]
\centering
\includegraphics[width=0.85\textwidth,trim={0cm 0cm 0cm 0cm},clip]{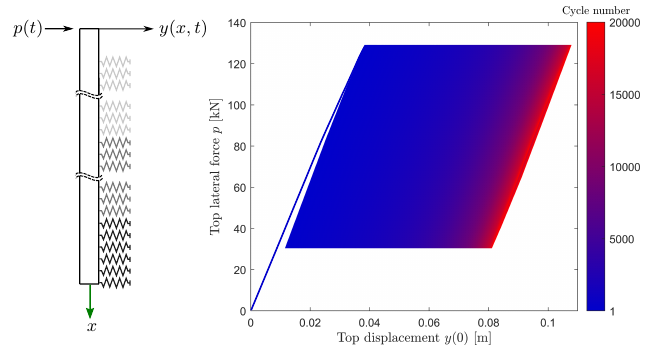}
\caption{Winkler pile foundation subjected to lateral cyclic loading (left) and corresponding force-displacement response (right). The pile is modeled as a circular hollow beam of 15 m depth, 0.92 m inner radius, and 1 m outer radius. The beam is connected to 45 equidistant springs that simulate the soil~support.}
\label{fig:pile1}
\end{figure}

The second example examines the behavior of a Winkler-beam monopile foundation under lateral cyclic loading (figure~\ref{fig:pile1}). The pile is modeled as an elastic Euler beam while the springs obey the ratcheting plasticity model discussed in section \ref{sec:const}, with a simplification to the one-dimensional case. The beam has a Young's modulus $E=210$~GPa and Poisson's ratio $\nu=0.3$. The springs are divided into 3 layers of 5 m depth with the following properties (in MPa): $E=266.67$, $\sigma^\mathrm{p}=2$, and $H^\mathrm{kin}=1466.7$ for the top layer ($0\leq x<5$~m); $E=1000$, $\sigma^\mathrm{p}=2.67$, and $H^\mathrm{kin}=2666.7$ for the middle layer ($5\leq x<10$~m); and $E=1333.3$, $\sigma^\mathrm{p}=3.33$, and $H^\mathrm{kin}=4666.7$ for the deepest layer ($10\leq x\leq15$~m). An isotropic hardening modulus $H^\mathrm{iso}=0$ and a ratcheting constant $\beta=0.01$ are used for all springs. 

The lateral load is applied for 20002 cycles of 101 load steps, with the load amplitude $p$ ranging from $30$~kN to $130$~kN. The force-displacement curve in figure~\ref{fig:plate2} shows a complex cyclic behavior, with non-constant ratcheting. A stable shakedown response is approached at the final cycles.

\begin{figure}[t!]
\centering
\includegraphics[width=1\textwidth,trim={0cm 0cm 0cm 0cm},clip]{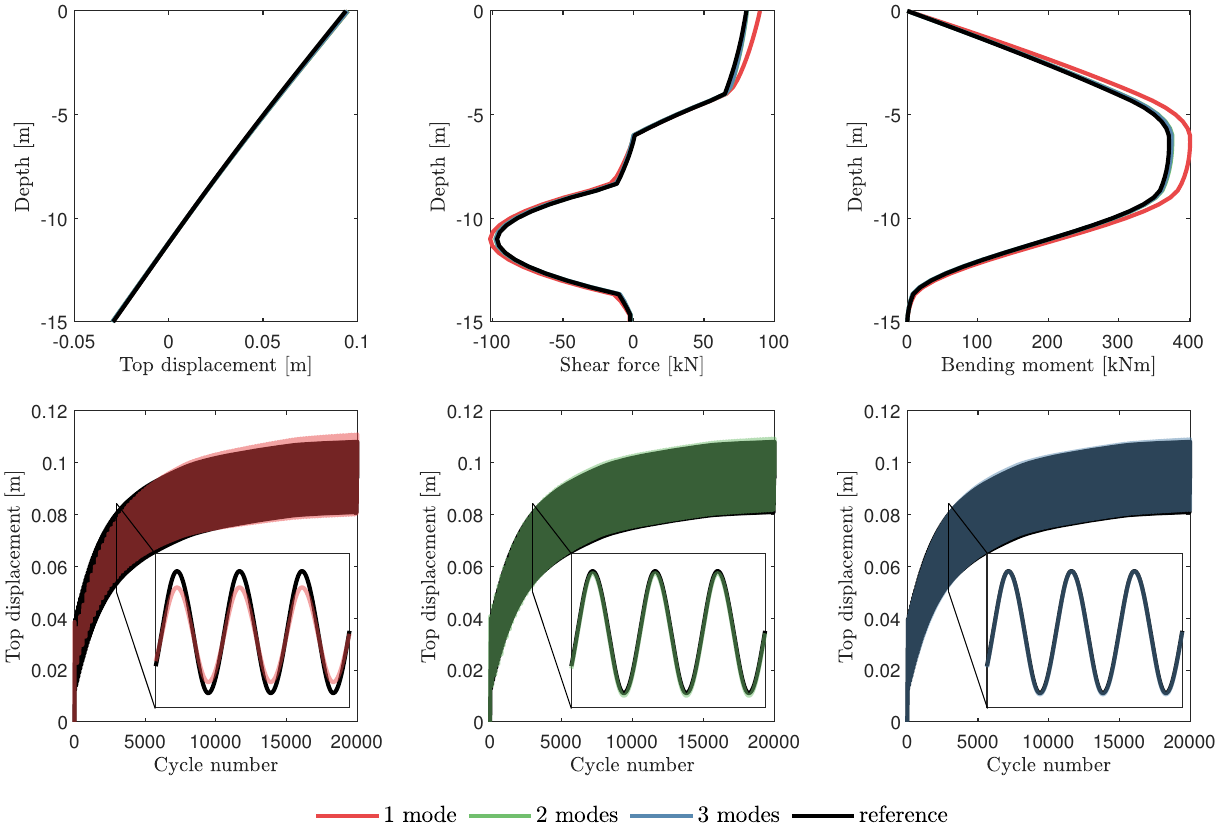}
\caption{Reference solution and PGD solution with 1, 2, and 3 modes: displacements, shear forces, and bending moments in the pile (top), and displacement time history~(bottom). PGD solution with 2 large-time scales ($N_1=200$ and $N_2=100)$.}
\label{fig:pile2}
\end{figure}

Similar to the previous example, we simulate two cycles in a standard incremental form, extrapolate the response periodically for 20000 cycles, and use the resulting time history as the initial guess for the PGD solution. We consider 3 temporal scales, with $N_\tau=101$ for the small-time scale and $N_1=200$ and $N_2=100$ for the 2 large-time scales. In this case, the number of temporal degrees of freedom is notably reduced. Assuming $M=3$, and noting that $N_\mathrm{d}=92$, we go from $1.84\times10^8$  space-time degrees of freedom in the classical incremental scheme to only $2.8776\times10^4$ in the PGD solution. 

Figure~\ref{fig:pile2} presents a comparison between the reference incremental response and the predicted PGD response using 1, 2, and 3 modes. The spatial distribution of displacements, shear forces, and bending moments at the final load cycle shows a relatively close agreement when considering a single mode and a very close agreement when considering 2 or 3 modes. The displacement time histories exhibit a similar trend. These results suggest that the proposed formulation holds promise as a numerical tool for studying monopile foundations, such as those found in offshore wind turbines, subjected to a large number of load cycles. However, it is important to note that the implemented algorithm, specifically the decoupled scheme in section~\ref{sec:numimp}, exhibited slow convergence from the initial guess (periodic extrapolation of the second cycle) to the final response. 

\FloatBarrier

\section{Conclusions and perspectives}

We have presented a multi-temporal PGD formulation for solids under cyclic loading, considering an arbitrary number of time scales. The proposed framework allows us to solve the spatial balance equations exclusively for the small-time intra-cyclic response. In contrast, the large-time inter-cyclic response is described by simple scalar algebraic equations. The framework is applied to a ratcheting plasticity model that exhibits strongly nonlinear cyclic behavior, incorporating competing mechanisms such as cyclic hardening and strain accumulation. The numerical results of relevant boundary value problems demonstrate that the PGD solutions accurately capture the reference response with a significantly reduced number of space-time degrees of freedom. Indeed, good approximations were achieved by employing only 3 modes in the examined problems. This approach holds promise for simulating complex structures and systems subjected to cyclic loads, such as offshore wind turbines. 

The present work has not addressed the pressing issue of computational time. Future research should focus on developing efficient numerical implementation schemes. \tB{For this purpose, a consistent linearization (section~\ref{sec:numimp}) may be considered, performing a space-time decomposition of the consistent tangent~\cite{arjoune2022} in a multi-temporal setting}. Additionally, efficient full-time integration schemes for the evolution equations of internal variables could be explored. For instance, the PGD formulation of the mechanical balance equations can be combined with a cycle jumping technique for the local evolution equations. Another approach to consider is a multi-temporal decomposition of the evolution equations, although the dissipative nature of internal variables may pose additional challenges. \tB{These potential developments would build upon the present study, enabling an efficient computation of non-linear responses spanning a substantial number of load cycles.} By doing so, we may fully harness the computational benefits offered by multi-temporal space-time PGD.

\tB{The proposed multi-temporal PGD also sets the stage for applications to structural systems under complex dynamic loading conditions. In this case, the loading functions may be separated into multiple time scales using a posteriori MOR techniques, such as the singular value decomposition. Then, a system similar to~\eqref{eq:pgd_spat_disc_lin0}--\eqref{eq:pgd_temp_disc_lin0} would arise but with ODEs replacing the algebraic equations describing the large-time~scales. This problem may be tackled by extending the space-time reformulation of Newmark's method \cite{boucinha2013,bamer2021} to the present multi-temporal case.}


\appendix
\FloatBarrier


\end{document}